# Comparison of continuity equation and Gaussian mixture model for long-term density propagation using semi-analytical methods

Pan Sun[1] · Camilla Colombo[2] · Mirko Trisolini[2] · Shuang Li[1]

**Abstract**
This paper compares the continuum evolution for density equation modelling and the Gaussian mixture model on the 2D phase space long-term density propagation problem in the context of high-altitude and high area-to-mass ratio satellite long-term propagation. The density evolution equation, a pure numerical and pointwise method for the density propagation, is formulated under the influence of solar radiation pressure and Earth's oblateness using semi-analytical methods. Different from the density evolution equation and Monte Carlo techniques, for the Gaussian mixture model, the analytical calculation of the density is accessible from the first two statistical moments (i.e. the mean and the covariance matrix) corresponding to each sub-Gaussian distribution for an initial Gaussian density distribution. An insight is given into the phase space long-term density propagation problem subject to nonlinear dynamics. The efficiency and validity of the density propagation are demonstrated and compared between the density evolution equation and the Gaussian mixture model with respect to standard Monte Carlo techniques.

## 1 Introduction

The possibility to exploit passive devices for the end-of-life disposal of high-altitude and high area-to-mass ratio satellites has been studied (Krivov and Getino 1997; Lücking et al. 2011a, b; Skoulidou et al. 2019). These satellites are usually small in size, restricted on the storage of on-board propellant and featured with little fuel left at the end of the mission, and are mainly perturbed by Solar Radiation Pressure (SRP) and Earth's oblateness ($J_2$) (Krivov and Getino 1997; Lücking et al. 2011a, b; Colombo and McInnes 2011a; Colombo et al. 2012; Skoulidou et al. 2019; Gkolias et al. 2020). As an example, the H2020 ReDSHIFT (Revolutionary Design of Spacecraft through Holistic Integration of Future Technologies) project boosted the research into passive end-of-life disposal, aiming at finding passive means to mitigate the proliferation of space debris ranging from Low Earth Orbit (LEO) to Geostationary Earth Orbit (GEO) (Rossi et al. 2018).

For the long-term dynamical evolution in the context of passive end-of-life disposal, the influence of orbit perturbations is essential to understand the motion of high-altitude and high area-to-mass ratio satellites and predict their disposal. To study the preliminary dynamical evolution, semi-analytical methods provide an efficient way to analyse the effect of orbit perturbations, since they filter out the short periodic terms of the disturbing function, retaining only long-term and secular terms. On the one hand, the numerical integration can be performed with a larger integration step, leading to a lighter computational effort. On the other hand, an insight into the dynamical characteristics is obtained as the periodic variation of the dynamics over one orbit revolution is averaged out. At the same time, for long-term propagation, an accuracy satisfying typical problem requirements is maintained (Colombo et al. 2012; Wittig et al. 2017; Rossi et al. 2018).

The problem of the long-term evolution of the density, both the probability density and uncertainty, was studied in several applications: the evolution of the interplanetary dust (Gor'kavyi et al. 1997a, b), the dynamics of nanosatellite constellations (McInnes 2000), high area-to-mass ratio spacecraft (Colombo and McInnes 2011b), the global debris population (Nazarenko 1997; Smirnov et al. 2001), and the evolution of debris cloud (Letizia et al. 2015; Letizia et al. 2016a, b; Frey 2020). The Density Evolution Equation (DEE) or continuity equation is a general method applied to the aforementioned long-term density propagation problems within different research domains. Its main idea is to consider the probability density as a fluid with continuous properties, which changes under the dynamics under consideration and can be obtained together with the propagation of the state variables. Thus, given the initial condition of the density distribution and the specific perturbation terms to be considered, the evolution of the density can be obtained with low computational effort. For the long-term phase space density propagation problem studied in this paper, Differential Algebra (DA) is another tool promising for the propagation of clouds of initial conditions in any sufficiently smooth dynamical system (Armellin et al. 2010; Wittig et al. 2017). Different from DEE, it can accomplish the density propagation with a single integration since it allows computation of a high-order polynomial expansion of the final density state as a function of the initial density state (Valli et al. 2013). In (Wittig et al. 2017), the authors introduce and combine for the first time the two techniques of semi-analytical methods and DA to allow efficient long-term evolution of a cloud of high area-to-mass ratio objects in Medium Earth Orbit (MEO).

The uncertainty propagation problem has been analysed in many research directions in astrodynamics, such as orbit determination (Lopez-Jimenez et al. 2020; Jia and Xin 2020), relative motion (Yang et al. 2018), planetary re-entry (Halder and Bhattacharya 2011; Jiang and Li 2019; Trisolini and Colombo 2021). A review of uncertainty propagation methods in orbit mechanics is given in Luo and Yang (2017). Apart from the traditional method of Monte Carlo (MC), which is used to provide reliable density propagation results with a large number of scattered samples, multiple nonlinear uncertainty propagation methods are available, such as Unscented Transformation (UT) (Julier et al. 1995, 2000; Julier and Uhlmann 2002, 2004), State Transition Tensor (STT) (Park and Scheeres 2006, 2007), Gaussian Mixture

Model (GMM) (Vittaldev et al. 2016; Frey 2020), Arbitrary Polynomial Chaos (APC) (Jia and Xin 2020), and the hybrid methods of Gaussian Mixture Model–Unscented Transformation (GMM-UT), Gaussian Mixture Model–State Transition Tensor (GMM-STT) (Yang et al. 2018), etc. Both UT and APC are sample-based uncertainty methods as MC, requiring the generation of a specific number of quasi-random samples. This is challenging when the uncertainty space dimension is high. The main idea of GMM is to approximate an arbitrary Probability Density Function (PDF) by a finite sum of weighted sub-Gaussian PDFs. Theoretically, the combined PDF can represent the real PDF with a large number of sub-Gaussians. For high-dimensional and nonlinear problems, as shown by Wittig et al. (2015), Jia and Xin (2020), Feng et al. (2021), the methodology of domain splitting can improve the accuracy of the uncertainty propagation.

The long-term density propagation problem has been dealt with in several aforementioned directions. However, to the authors' knowledge, few studies have been published on the long-term phase space density propagation problem. Wittig et al. 2017 focus on the differential algebra and semi-analytical method for the long-term phase space density propagation problem, but the application of a pure numerical and pointwise method retaining nonlinear characteristics during density propagation is not presented. Different from their research works, this paper compares density evolution equation (featuring density calculation from a pure numerical view) and Gaussian mixture model (allowing analytical calculation of the density) with respect to MC on the 2D phase space long-term density propagation problem. The density evolution equation is formulated under the influence of solar radiation pressure and Earth's oblateness using the semi-analytical method (Colombo and McInnes 2011b). The technique of linear interpolation is introduced for accurate density calculation with DEE. Unscented transformation is used to nonlinearly propagate the first two statistical moments with GMM.

The paper is organised as follows: Sect. 2 gives the problem formulation and analysis, including the semi-analytical equation (focusing on Hamiltonian phase space) and the formulation of the density evolution equation. Section 3 presents detailed method description, and the density propagation methods and techniques of MC, DEE and GMM-UT, respectively. Section 4 presents the simulation setup. In Sect. 5, density evolution results are presented for the 2D phase space. Discussion is given on the comparison of the density evolution results within nonlinear phase space domains and on the comparison of the computational accuracy and efficiency of DEE and GMM-UT with respect to MC. Section 6 presents some conclusions.

## 2 Problem formulation and analysis

### 2.1 Semi-analytical equation

To study the dynamical evolution under the influence of SRP and $J_2$, a semi-analytical method is used in this paper (Krivov and Getino 1997; Lücking et al. 2011a, 2011b).

The semi-analytical equations of motion for the eccentricity $e$ and the solar angle $\phi$ (see detailed geometrical representation of the problem in "Appendix A"), initially derived in studies of circumplanetary dust dynamics in (Hamilton 1993) and (Hamilton and Krivov

1996), are as follows:

$$\begin{cases} \dfrac{de}{dt} = n_{sun}\left(C\sqrt{1-e^2}\sin\phi\right) \\ \dfrac{d\phi}{dt} = n_{sun}\left(C\dfrac{\sqrt{1-e^2}}{e}\cos\phi + \dfrac{W}{(1-e^2)^2} - 1\right), \end{cases} \quad (1)$$

where $C$ and $W$ represent the dimensionless radiative and oblateness parameters, respectively, and $n_{sun}$ is the mean motion of the sun. The definitions of $C$ and $W$ are (Hamilton and Krivov 1996; Colombo et al. 2012)

$$\begin{cases} C = \dfrac{3}{2}\sigma\dfrac{n_s}{n_{sun}} \\ W = \dfrac{3}{2}J_2\left(\dfrac{R_E}{a}\right)^2\dfrac{n_s}{n_{sun}}, \end{cases} \quad (2)$$

where $n_s$ is the mean motion of the satellite, $J_2$ and $R_E$ are the second-order zonal harmonic coefficient and the equatorial radius of the Earth, respectively, and $\sigma$ is the ratio of radiative force to Earth's gravity, $\sigma = F_\odot a^2 \tau/(\mu c)$, where $F_\odot$ is the solar flux at 1 AU, $\mu$ is the gravitational parameter of the Earth, $c$ is the speed of light in vacuum, and $\tau$ is the area-to-mass ratio of the satellite.

The quasi-canonical form of Eq. (1) can be written as follows (Hamilton 1993; Hamilton and Krivov 1996),

$$\begin{cases} \dfrac{de}{dt} = n_{sun}\left(-\dfrac{\sqrt{1-e^2}}{e}\dfrac{\partial H}{\partial \phi}\right) \\ \dfrac{d\phi}{dt} = n_{sun}\left(\dfrac{\sqrt{1-e^2}}{e}\dfrac{\partial H}{\partial e}\right), \end{cases} \quad (3)$$

where $H$ is the phase space Hamiltonian, $H = \sqrt{1-e^2} + Ce\cos\phi + \dfrac{W}{3}(1-e^2)^{-3/2}$. $H(\phi, e)$ = constant, i.e. for the phase space variable pair $(\phi, e)$ defined at time $t$, the phase space evolution is along a constant Hamiltonian contour line.

To determine the phase space structure and evolution properties, the stationary points of the Hamiltonian $H$ are calculated by solving for $\partial H/\partial \phi = \partial H/\partial e = 0$. Summarised results of stationary points constrained to the planar problem are given in (Krivov and Getino 1997), while solutions extended to non-planar orbits can be found by referring to the work of Gkolias et al. 2020.

With the planar semi-analytical equations and solutions of stationary points, particular attention is paid to the application of the passive end-of-life disposal in the works of Krivov and Getino 1997 and Lücking et al. 2011a. In Krivov and Getino 1997, a detailed classification of the phase portraits is given. Without imposing restrictions on eccentricity and area-to-mass ratio, in (Krivov and Getino 1997) the authors present the evolution of the eccentricity and apse line for a satellite with varying area-to-mass ratios and geocentric distances. Figure 1 shows one phase portrait of type III (Krivov and Getino 1997), including five stationary points, obtained with the initial conditions of $a = 2.5\,R_E$, $W = 0.409$, $C = 0.15$, similar to Fig. 2f in (Krivov and Getino 1997). The horizontal line marks the critical eccentricity $e_{cri} = 0.6$ for Earth re-entry at the Earth surface. In this phase portrait, phase space bifurcation is detected at the stationary point $P_4$. The Hamiltonian phase space is divided into three sub-domains departed by the contour lines passing by the stationary points $P_1$ and $P_4$. As shown in Fig. 1,

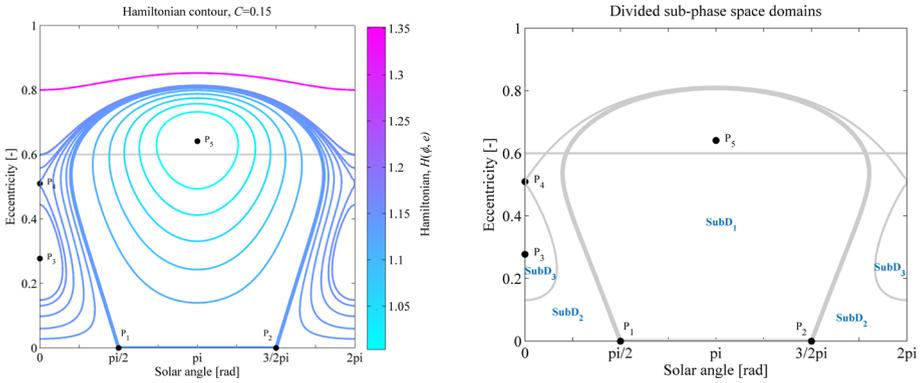

**Fig. 1** Phase portrait of type III; Divided sub-phase space domains ($a = 2.5\,R_E$, $W = 0.409$, $C = 0.15$, critical eccentricity $e_{cri} = 0.6$ for Earth re-entry)

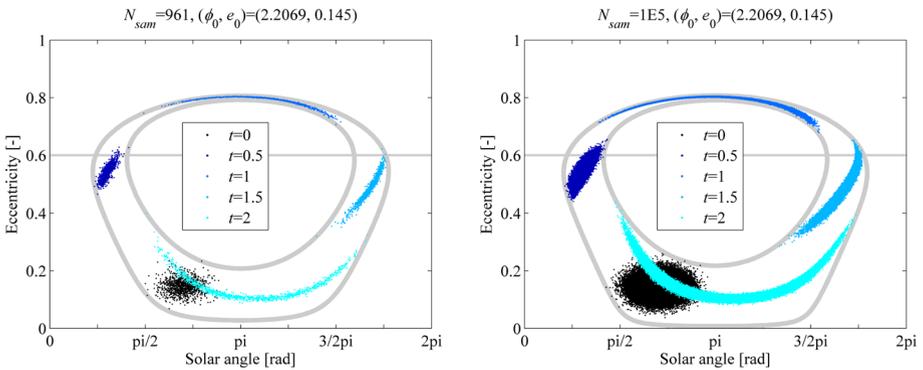

**Fig. 2** Phase space evolution for Scenario 1 within $SubD_1$, $N_{sam} = 961$ vs. $N_{sam} = 1E5$ (critical eccentricity $e_{cri} = 0.6$ for Earth re-entry)

we define the three sub-domains as $SubD_i$, $i \in \{1,2,3\}$, residing between the Hamiltonian contour lines $H_{P_5}$ and $H_{P_1}$, $H_{P_1}$ and $H_{P_4}$, $H_{P_4}$ and $H_{P_3}$, respectively. Since the dynamical characteristics change quickly and greatly near the bifurcation point, it is of great interest to understand how the phase space evolves with time when density information is considered (i.e. understand how the phase space and its associated density evolve with time) within each sub-domain. It should be noted that only the intermediate zone ($2R_E \leq a \leq 3R_E$) (Krivov and Getino 1997) best applies the semi-analytical model under the influence of SRP and $J_2$. For satellites with $a \geq 3\,R_E$ (i.e., with altitudes larger than 13,000 km), the model using only SRP and $J_2$ is no longer adequate due to the significant lunisolar perturbations. In this case, the coupling SRP and lunisolar perturbation problem should be examined (Hamilton and Krivov 1996). For satellite orbits with both high apogees and low perigees, the effects of the solar apsidal resonance (due to the interplays of the lunisolar perturbation, Earth's oblateness, and air drag) on orbital evolution should be considered (Wang and Gurfil 2017). In this paper, we analyse the evolution of the phase space density to study the reliability and robustness of passive end-of-life disposal solutions for high-altitude and high area-to-mass ratio satellites within the intermediate zone ($2R_E \leq a \leq 3R_E$) under the influence of SRP and $J_2$ only. The focus of this paper is on the comparison of the density evolution equation and Gaussian mixture model for long-term density propagation. It is true that adding the effect of lunisolar perturbation would improve the accuracy of the results, but this would not change the outcome of the work which aims at comparing different uncertainty propagation techniques applied to a semi-analytical dynamics. The coupling with the solar gravitational perturbation that can occur at higher altitudes (and for a low perigee, the coupling terms from the air drag) will be given into an insight in the future work.

### 2.2 Density evolution equation

To study the density evolution in the phase space, we rely on a continuum method, which propagates the density evolution equation together with the semi-analytical equations of motion.

Assuming that $n$ is the density to be solved for a specific problem, the density evolution equation (Gor'kavyi et al. 1997a; McInnes 2000; Frey 2020) can be written as,

$$\frac{\partial n}{\partial t} + \nabla \cdot \mathbf{f} = \dot{n}^+ - \dot{n}^-, \tag{4}$$

where $\nabla \cdot \mathbf{f}$ represents the continuous acceleration terms to be considered for the dynamical system, such as the perturbation of Earth's oblateness in this paper, and $\dot{n}^+ - \dot{n}^-$ represents the discontinuous terms to be considered, such as the injection of new fragments due to launches on the topic of debris population evolution.

Given $m$ generic variables $\alpha_i, i \in \{1, \cdots, m\}$, and assuming that the density $n$ is differentiable with respect to all $\alpha_i$ everywhere (Letizia et al. 2016a, b), Eq. (4) can be rewritten in generic coordinates as follows.

$$\frac{\partial n}{\partial t} + \frac{\partial n}{\partial \alpha_1} v_{\alpha_1} + \cdots + \frac{\partial n}{\partial \alpha_m} v_{\alpha_m} + \left[ \frac{\partial v_{\alpha_1}}{\partial \alpha_1} + \cdots + \frac{\partial v_{\alpha_m}}{\partial \alpha_m} \right] n = \dot{n}^+ - \dot{n}^-. \tag{5}$$

In this paper, $\dot{n}^+ - \dot{n}^- = 0$. In this way, by applying the method of characteristics (Evans 1998), the following Ordinary Differential Equations (ODEs) are obtained (Letizia et al. 2016a, b),

$$\begin{cases} \dfrac{dt}{ds} = 1 \\ \dfrac{d\alpha_1}{ds} = v_{\alpha_1}(\alpha_1, \cdots, \alpha_m) \\ \vdots \\ \dfrac{d\alpha_m}{ds} = v_{\alpha_m}(\alpha_1, \cdots, \alpha_m) \\ \dfrac{dn}{ds} = -\left[ \dfrac{\partial v_{\alpha_1}}{\partial \alpha_1} + \ldots + \dfrac{\partial v_{\alpha_m}}{\partial \alpha_m} \right] n(\alpha_1, \ldots, \alpha_m, t), \end{cases} \tag{6}$$

where $s$ represents a parameterisation of the characteristic lines. It can be seen that given the specific formulation of the problem, i.e. given the actual expressions of $v_{\alpha_i}$, the result of the density $n(\alpha_1, \cdots, \alpha_m, t)$ can be obtained by numerical integration together with the state variables $\alpha_i$. Apart from the numerical way to solve the density evolution equation (Trisolini and Colombo 2021), analytical solutions can be found for specific problems (McInnes 2000; Colombo and McInnes 2011b; Letizia et al. 2015).

For the case in examination, to avoid singularities during numerical integration, we rewrite the equations of motion of Eq. (3) as a function of two new variables $x_1 = e \cdot \sin(\phi)$ and $x_2 = e \cdot \cos(\phi)$ as follows:

$$\begin{cases} \dfrac{dx_1}{dt} = \dfrac{de}{dt}\sin(\phi) + e\cos(\phi)\dfrac{d\phi}{dt} \\ \dfrac{dx_2}{dt} = \dfrac{de}{dt}\cos(\phi) - e\sin(\phi)\dfrac{d\phi}{dt}. \end{cases} \quad (7)$$

Substituting Eq. (1) into Eq. (7) and expressing all terms concerned with $e$ and $\phi$ in $x_1$ and $x_2$, i.e., $e = (x_1 + x_2)^{1/2}$, $\sin(\phi) = x_1/e$, $\cos(\phi) = x_2/e$, the expressions of $v_{\alpha_i}$ in Eq. (6) are derived. In this way, the differential equation controlling the density evolution becomes,

$$\frac{dn}{dt} = -\left(\frac{\partial}{\partial x_1}\left(\frac{dx_1}{dt}\right) + \frac{\partial}{\partial x_2}\left(\frac{dx_2}{dt}\right)\right)n = \frac{n_{sun} \cdot x_1 \cdot C}{\left(-x_1^2 - x_2^2 + 1\right)^{1/2}} n. \quad (8)$$

From Eq. (8), we can see that the density evolution under the influence of solar radiation pressure and Earth's oblateness is independent of the oblateness parameter $W$, i.e. the existence of the influence of Earth's oblateness alone does not change the density with time, which is consistent with the analysed result of the case of Earth's oblateness only in (Letizia et al. 2016a, b).

## 3 Density propagation methods and techniques

Given an initial state and the corresponding density distribution, the main idea of density evolution is to obtain the density distribution corresponding to the final state. Both MC and DEE give an insight into intrinsic nonlinear characteristics of the dynamical system during the density evolution. Different from MC and DEE, for GMM-UT, analytical calculation of the density can be achieved with the propagated first two statistical moments corresponding to each sub-Gaussian. To analyse different methods, detailed method description and density calculation techniques are given as follows for MC, DEE and GMM-UT, respectively.

### 3.1 Monte Carlo

In this paper, the MC method is used as a reference to validate and compare the performance of DEE and GMM-UT. For MC, it is enforced to generate and propagate a large number of random samples. The number of samples varies depending on multiple aspects, such as the nonlinearity of the dynamics, the dimension of the studied state space and the specified confidence interval (Wilson 1927; Jehn 2015; Wallace 2015; Letizia et al. 2016a, b; Romano 2020). To obtain the joint and marginal density for the case in examination for MC, a binning approach is used.

#### 3.1.1 Binning approach

Assume a generic 2D problem defined in two independent variables $x_i$, $i \in \{1,2\}$, and the total number of the propagated samples at a specific time instant $t$ is $N_{sam}$. The main idea of the binning approach is to partition values in $x_1$ and $x_2$ into the 2D uniformly divided bins, with defined number of bins, $N_{bi}$, and the bin edges in each dimension, $Edges_i$. The minimum and maximum values of the bin edges, $LE_i$ and $UE_i$, of $Edges_i$, are determined by the minimum and maximum value of $x_i$. Thus, the bin width for each dimension is $wid_i = (UE_i\text{-}LE_i)/N_{bi}$, and the element set of $Edges_i$ is determined, i.e. $Edges_i(ct) = LE_i + wid_i \cdot (ct\text{-}1)$, $ct \in \{1,\cdots, N_{bi}+1\}$.

### 3.1.2 Joint and marginal density calculation

With the defined 2D bins, the number of samples in each bin $c_{pk}$, $p \in \{1,\cdots, N_{b1}\}$, $k \in \{1,\cdots, N_{b2}\}$, is counted. Then, the 2D joint density of MC is obtained as follows,

$$f_{MC-pk}(x_1, x_2, t) = c_{pk}/(N_{sam} \cdot A_{pk}), \tag{9}$$

where $A_{pk}$ is the area of each bin. As shown in Eq. (9), the joint density of MC represents the frequency of samples per bin area.

The marginal density for one dimension can be obtained by integrating the joint density throughout the whole domain of the other dimension. Here the calculation formulation of the marginal density with respect to the first dimension is presented.

$$f_{MC-1-p}(x_1, t) = \frac{A_{pk} \cdot \sum_{x_2=LE_2}^{UE_2} f_{MC-pk}(x_1, x_2, t)}{wid_1}, \quad p = \{1, \ldots, N_{b1}\}. \tag{10}$$

As indicated in Eq. (10), the marginal density in terms of $x_i$ of MC represents the frequency of samples per bin width of $wid_i$, $i \in \{1,2\}$.

## 3.2 Density evolution equation

Four steps are required to obtain the density with DEE. First, generate initial random samples in the 2D phase space subject to the defined initial density distribution with a predefined sample number. Second, calculate the initial probability density weights corresponding to the samples. Thus, initial samples given in the 2D phase space ($\phi$, $e$) and their associated density weights are obtained in the 3D state space ($\phi$, $e$, $n$). Third, integrate the density evolution equation together with the semi-analytical equations of motion to obtain final samples and their corresponding density weights in the 3D state space. Fourth, calculate the density by processing the final samples and density weights in a statistical way. In this paper, the technique of linear interpolation is used to calculate the density together with the binning approach.

### 3.2.1 Linear interpolation

For the specific 2D problem, the linear interpolation method based on Delaunay triangulation (Preparata and Shamos 1985) is adopted. The main advantage of this method lies in the capability of retaining the scattered sample data at the nodes of the triangulation. From this, it can be preliminarily concluded that with a larger number of scattered sample data, the quality of the linearly interpolated results may be better. Note that for the linear interpolation method based on Delaunay triangulation, the interpolation is done within the whole convex hull of the scattered sample data.

Give the sample data and their corresponding density weights at a specific time instant $t$ in the vectors $(x_{10}(t), x_{20}(t), n_0(x_{10}, x_{20}, t))$ and assume the grid number of $Ngrid_i, i \in \{1,2\}$, for each dimension. The fundamental of linear interpolation is summarised with the following mapping relationship.

$$n(x_1, x_2, t) = \text{Linear Interpolation}(x_{10}, x_{20}, n_0, Ngrid_1, Ngrid_2), \quad (11)$$

where $(x_1, x_2)$ are the query points specified to return the linearly interpolated results.

### 3.2.2 Binning approach

With the interpolated results in the vectors $(x_1(t), x_2(t), n(x_1, x_2, t))$ at the specific time instant $t$, the binning approach same as that presented for MC is introduced. However, for DEE, the binning approach is utilised to partition values in the density weights $n(x_1, x_2, t)$ into the 2D uniformly divided bins in the vectors $(x_1(t), x_2(t))$.

### 3.2.3 Joint and marginal density calculation

With the defined 2D bins, the mean value of the density weights $n(x_1, x_2, t)$ in each bin $N_{pk}$, $p \in \{1, \cdots, N_{b1}\}, k \in \{1, \cdots, N_{b2}\}$, is calculated. Then, the 2D joint density of DEE is obtained as follows,

$$f_{DEE-pk}(x_1, x_2, t) = \frac{N_{pk}/A_{pk}}{\sum_{p=1}^{N_{b1}} \sum_{k=1}^{N_{b2}} N_{pk}}, \quad (12)$$

where $A_{pk}$ is the area of each bin, and $N_{bi}, i \in \{1, 2\}$, is the number of bins for each dimension. As shown in Eq. (12), the joint density of DEE represents the weighted mean of the density weights per bin area.

Similar to MC, the marginal density for one dimension can be obtained by integrating the joint density throughout the whole domain of the other dimension. Here the calculation formulation of the marginal density with respect to the first dimension is presented.

$$f_{DEE-1-p}(x_1, t) = \frac{A_{pk} \cdot \sum_{x_2=LE_2}^{UE_2} f_{DEE-pk}(x_1, x_2, t)}{wid_1}, \ p = \{1, \ldots, N_{b1}\}. \quad (13)$$

As indicated in Eq. (13), the marginal density in terms of $x_i$ of DEE represents the weighted mean of the density weights per bin width of $wid_i, i \in \{1, 2\}$.

### 3.3 Gaussian mixture model-unscented transformation

For GMM-UT, four steps are required to obtain the density. First, split an initial Gaussian density distribution into $N_{1d}$ sub-Gaussians. Second, obtain initial sigma points of UT and propagate them in the 2D phase space. Third, calculate the transformed first two statistical moments corresponding to each sub-Gaussian with the propagated sigma points using UT. Then, get the combined joint and marginal density with the weighted sum of the propagated first two statistical moments corresponding to each sub-Gaussian.

To obtain the joint and marginal density of GMM-UT, the techniques of Gaussian mixture model and unscented transformation are introduced.

### 3.3.1 Gaussian mixture model

A Gaussian mixture model is a weighted sum of Gaussian PDFs (DeMars et al. 2013; Vittaldev et al. 2016; Vittaldev and Russell 2016; Sun and Kumar 2016; Yang et al. 2018),

$$f_c(\boldsymbol{x}, t; \boldsymbol{m}_c, \boldsymbol{P}_c) = \sum_{i=1}^{N_{1d}} \omega_i f_{gi}(\boldsymbol{x}, t; \boldsymbol{m}_{si}, \boldsymbol{P}_{si}), \qquad (14)$$

where $f_c(\boldsymbol{x}, t; \boldsymbol{m}_c, \boldsymbol{P}_c)$ is the combined PDF, $\boldsymbol{x}$ is the state vector, $\boldsymbol{m}_c$ and $\boldsymbol{P}_c$ are the combined mean vector and covariance matrix, respectively, $f_{gi}(\boldsymbol{x}, t; \boldsymbol{m}_{si}, \boldsymbol{P}_{si})$, $\boldsymbol{m}_{si}$, $\boldsymbol{P}_{si}$ and $\omega_i$ are the PDF, mean vector, covariance matrix and weight corresponding to the $i$th sub-Gaussian, respectively, and $\sum_{i=1}^{N_{1d}} \omega_i = 1$, $\omega_i \in (0,1]$. The weight $\omega_i$ can be calculated by minimising the difference between the combined PDF $f_c(\boldsymbol{x}, t; \boldsymbol{m}_c, \boldsymbol{P}_c)$ and the initial Gaussian PDF $f_t(\boldsymbol{x}, t; \boldsymbol{m}_t, \boldsymbol{P}_t)$. Theoretically, with a large number of sub-Gaussians (in the $L_1$ norm sense) (Alspach and Sorenson 1972), the combined GMM can approximate the initial Gaussian PDF. Since a Gaussian distribution can be completely determined by the mean and covariance matrix, only the first two statistical moments need to be propagated. Multiple nonlinear methods are available for nonlinear propagation of the first two statistical moments, such as UT, STT and Fourier–Hermite series (Mikulevicius and Rozovskii 2000; Sarmavuori and Sarkka 2011).

Before splitting an initial Gaussian density distribution, an univariate GMM library of the 1D standard Gaussian distribution is formed (Huber et al. 2008; Horwood et al. 2011; DeMars et al. 2013; Vittaldev and Russell 2013; Vittaldev et al. 2016), which is computed only once and is stored in terms of the means, weights and standard deviations. Univariate Gaussian splitting libraries splitted within three rules for odd number of splits up to 39 are provided in (Vittaldev and Russell 2016), in terms of the parameter set of the splitting number, means, weights and standard deviations, $\{N_{1d}, m_{i1d}, \omega_{i1d}, \sigma\}, i \in \{1, \cdots, N_{1d}\}$. Note that the standard deviation for each sub-Gaussian is the same here.

Given the initial Gaussian distribution $f_t(\boldsymbol{x}, t; \boldsymbol{m}_t, \boldsymbol{P}_t)$ subject to a specific dynamical problem, to apply the univariate Gaussian splitting library in (Vittaldev and Russell 2016), a splitting direction is determined along the most nonlinear one of the dynamical system. Assume that the $j$th direction of the state vector $\boldsymbol{x}$ is selected for splitting. For the $i$th sub-Gaussian, $i \in \{1, \cdots, N_{1d}\}$, the mean $\boldsymbol{m}_{si}$ and covariance matrix $\boldsymbol{P}_{si}$ are calculated as follows (DeMars et al. 2013),

$$\begin{cases} \omega_i = \omega_{i1d}, \; \boldsymbol{m}_{si} = \boldsymbol{m}_t + \sqrt{\lambda_j} m_{i1d} \boldsymbol{v}_j, \; \boldsymbol{P}_{si} = \boldsymbol{V} \boldsymbol{\Lambda}_i \boldsymbol{V}^T \\ \boldsymbol{\Lambda}_i = \mathrm{diag}\{\lambda_1, \lambda_2, \ldots, \sigma^2 \lambda_j, \ldots, \lambda_{Nvar}\}, \; \boldsymbol{P}_t = \boldsymbol{V} \boldsymbol{\Lambda} \boldsymbol{V}^T \end{cases}, \qquad (15)$$

where $\boldsymbol{\Lambda}_i$ is a diagonal matrix made up of the diagonal elements of $\boldsymbol{\Lambda}$ and the standard deviation $\sigma$ of univariate splitting, $\boldsymbol{\Lambda} = \mathrm{diag}\{\lambda_1, \lambda_2, \cdots, \lambda_{Nvar}\}$, $\boldsymbol{V} = [\boldsymbol{v}_1, \boldsymbol{v}_2, \cdots, \boldsymbol{v}_{Nvar}]$, and $Nvar$ is the state space dimension. Note that the only difference of $\boldsymbol{\Lambda}_i$ from $\boldsymbol{\Lambda}$ is the content of the $j$th diagonal element. $\boldsymbol{\Lambda}$ and $\boldsymbol{V}$ are the eigenvalue matrix and the eigenvector matrix of $\boldsymbol{P}_t$, respectively.

With the obtained $N_{1d}$ sub-Gaussians, $\{\omega_i, \boldsymbol{m}_{si}, \boldsymbol{P}_{si}\}, i \in \{1, \cdots, N_{1d}\}$, the nonlinear methods such as UT, STT and Fourier–Hermite series are used to propagate the mean and covariance matrix, $\{\boldsymbol{m}_{si}, \boldsymbol{P}_{si}\}$, of each sub-Gaussian. The combined non-Gaussian PDF, $f_c(\boldsymbol{x}, t; \boldsymbol{m}_c, \boldsymbol{P}_c)$, is calculated with Eq. (14), and the combined mean and covariance matrix, $\{\boldsymbol{m}_c, \boldsymbol{P}_c\}$, are calculated as follows (DeMars et al. 2013).

$$\begin{cases} \boldsymbol{m}_c = \sum_{i=1}^{N_{1d}} \frac{\omega_i}{\omega_m} \boldsymbol{m}_{si}, \ \omega_m = \sum_{i=1}^{N_{1d}} \omega_i \\ \boldsymbol{P}_c = \sum_{i=1}^{N_{1d}} \frac{\omega_i}{\omega_m} \left( \boldsymbol{P}_{si} + \boldsymbol{m}_{si} \boldsymbol{m}_{si}^T \right) - \boldsymbol{m}_c \boldsymbol{m}_c^T. \end{cases} \quad (16)$$

In this paper, the weights are set to be constant in the long-term density propagation, and the splitting is done only at the initial time along the more nonlinear direction of the solar angle. The choice of splitting in the solar angle direction only is determined by the nonlinearity measure formula of Eq. (5) given in (Vittaldev and Russell 2016). The technique of unscented transformation is used to nonlinearly propagate the first two statistical moments of each sub-Gaussian.

### 3.3.2 Unscented transformation

The main idea of unscented transformation is that approximating the probability distribution is less complicated than approximating the nonlinear transformation (Julier et al. 1995), which means by deterministically choosing and integrating a few samples, the probability distribution at a specific time can be approximated with the transformed probability moments calculated with the few propagated samples.

In this paper, the symmetric extended set (Julier and Uhlmann 2004) is employed for the determination of the initial $N_s$ sigma points, where $N_s = (2 \cdot Nvar + 1)$. Three steps are required to calculate the transformed first two statistical moments corresponding to each sub-Gaussian. First, the set of initial weighted sigma points, $\boldsymbol{x}_{sk}(t_0), k \in \{0, 1, 2, \cdots, 2 \cdot Nvar\}$, is deterministically chosen such that the true mean $\boldsymbol{m}_s(t_0)$ and covariance matrix $\boldsymbol{P}_s(t_0)$ at the initial time $t_0$ for each sub-Gaussian are precisely captured. Here the initial $N_s$ sigma points are defined as follows,

$$\begin{cases} \boldsymbol{x}_{s0} = \boldsymbol{m}_s(t_0) \\ \boldsymbol{x}_{si} = \boldsymbol{m}_s(t_0) + \boldsymbol{S}_{si} \cdot \sqrt{Nvar+\varsigma} \\ \boldsymbol{x}_{s(2i)} = \boldsymbol{m}_s(t_0) - \boldsymbol{S}_{s(2i)} \cdot \sqrt{Nvar+\varsigma}, \ i \in \{1, 2, \ldots, Nvar\}, \end{cases} \quad (17)$$

where $\boldsymbol{S}_s$ is the square root matrix of the covariance matrix $\boldsymbol{P}_s(t_0)$, i.e., $\boldsymbol{P}_s(t_0) = \boldsymbol{S}_s \boldsymbol{S}_s^T$, $\boldsymbol{S}_{si}$ is the $i$th column of $\boldsymbol{S}_s$, $\zeta = \alpha^2(Nvar + \beta) - Nvar$, and $\alpha$ and $\beta$ are two parameters that determine the distribution of sigma points around the mean value. In this paper, $\alpha = 0.8$ and $\beta = 0$ are assumed. Second, propagate the initial selected sigma points through the nonlinear dynamical equations to obtain the nonlinearly transformed sigma points at a specific time instant $t$, $\boldsymbol{x}_{sk}(t), k \in \{0, 1, 2, \cdots, 2 \cdot Nvar\}$. Then, the transformed mean and covariance matrix are calculated as follows (Julier and Uhlmann 2004),

$$\begin{cases} \boldsymbol{m}_s(t) = \sum_{k=0}^{N_s-1} w_k^m \boldsymbol{x}_{sk}(t) \\ \boldsymbol{P}_s(t) = \sum_{k=0}^{N_s-1} w_k^P [\boldsymbol{x}_{sk}(t) - \boldsymbol{m}_s(t)][\boldsymbol{x}_{sk}(t) - \boldsymbol{m}_s(t)]^T, \end{cases} \quad (18)$$

where $w^m_k$ and $w^P_k$ are the weights for calculating the transformed mean and covariance matrix with the $k$th propagated sigma point, respectively, $k \in \{0, 1, 2, \cdots, 2 \cdot Nvar\}$. The specific weights for the $N_s$ sigma points are given as follows,

$$\begin{cases} w^m_0 = \dfrac{\varsigma}{Nvar+\varsigma}, k=0 \\ w^P_0 = \dfrac{\varsigma}{Nvar+\varsigma} + (1-\alpha^2+\eta), k=0 \\ w^m_k = w^P_k = \dfrac{\varsigma}{2(Nvar+\varsigma)}, k=\{1,2,\ldots,2\cdot Nvar\}, \end{cases} \quad (19)$$

where $\eta$ is a scale parameter introduced to determine the distribution of sigma points around the covariance matrix. In this paper, $\eta=2$ is assumed.

### 3.3.3 Joint and marginal density calculation

For the case in examination, given an initial Gaussian distribution $f_t(\mathbf{x}, t; \mathbf{m}_t, \mathbf{P}_t)$, with the initially splitted $N_{1d}$ sub-Gaussians given by Eq. (15) and the propagated first two statistical moments of each sub-Gaussian given by Eq. (18) using UT, the combined joint density $f_{GMM\text{-}UT}(x_1, x_2, t)$ is determined with Eq. (14) by calculating the weighted sum of the joint density corresponding to each sub-Gaussian.

Since the first two statistical moments corresponding to each sub-Gaussian are given by Eq. (18), the marginal density is determined by calculating the weighted sum of the 1D Gaussian PDFs using the means and standard deviations corresponding to sub-Gaussians along the desired direction.

## 4 Simulation setup

### 4.1 Definition of initial conditions and case content

To understand how the phase space and its associated density evolve with time in each sub-phase space domain shown in Fig. 1, three scenarios are considered (Table 1). Each scenario defines a different initial condition in terms of a Gaussian density distribution. $\phi_0$ and $e_0$ represent the mean of the Gaussian distribution in the $(\phi, e)$ phase space, while $\Delta \phi$ and $\Delta e$ specify the covariance matrix $P_0 = [(\Delta \phi/2)^2\ 0; 0\ (\Delta e/2)^2]$, to make sure of the realisation of the long-term density propagation within each sub-phase space domain. $T_U$ is the simulation time and $dt$ represents the time interval between two consecutive snapshots.

For each scenario, four cases are considered to validate the feasibility of the formulated density evolution equation and compare the performance of DEE and GMM-UT with MC. A detailed description for the four cases is given in Table 2, where $N_{sam}$ and $N_{1d}$ are the number of initial samples for MC and DEE, and the splitting number of the initial Gaussian density distribution for GMM-UT, respectively. *Ngrid* defines the grid number to perform linear interpolation in two phase space dimensions for both Case 2 and Case 3 for DEE. The number of MC random samples is set to be 1E5. The same sample number $N_{sam}$=1E5 as that of MC and a larger grid number *Ngrid*=5E3 for Case 3 are set to compare the accuracy of DEE with MC. The much smaller sample number $N_{sam}$=961 and much smaller gird number *Ngrid*=1E3 of Case 2 compared with that of Case 3 are set to compare the accuracy and efficiency of DEE with MC and GMM-UT. The smaller sample number $N_{sam}$=961 is determined to make sure of the accuracy of the density calculated using the linear interpolation method, accounting for about 1% of the larger sample number $N_{sam}$=1E5. The splitting number of GMM-UT is set to be 39 using the univariate splitting library of the maximum splitting number provided by (Vittaldev and Russell 2016), and the more nonlinear direction of the solar angle is set as the splitting direction.

**Table 1** Definition of initial conditions for the three test case scenarios ($a = 2.5\ R_E$)

| Scenario | $\phi_0$, rad | $e_0$ | $\Delta \phi$, rad | $\Delta e$ | $T_U$, yr | $dt$, yr |
|---|---|---|---|---|---|---|
| 1 | 2.2069 | 0.145 | $\pi/8$ | 0.05 | 2 | 0.5 |
| 2 | 0.5419 | 0.095 | $\pi/40$ | 0.01 | 3 | 0.5 |
| 3 | 0.3004 | 0.23 | $\pi/32$ | 0.02 | 2 | 0.5 |

**Table 2** Definition of the four cases for each scenario

| Case | 1 | 2 | 3 | 4 |
|---|---|---|---|---|
| Identifier | MC | DEE-961 | DEE-1E5 | GMM-UT |
| $N_{sam}$ | 1E5 | 961 | 1E5 | – |
| $N_{1d}$ | – | – | – | 39 |
| Ngrid | – | 1E3 | 5E3 | – |

### 4.2 Definition of computational effort

To compare the computational effort of DEE and GMM-UT with MC, the main contribution of the computational effort is analysed for the three methods and the main parts of the computational effort are defined.

The main computational effort of MC includes two aspects. First, generate initial random samples at $t = 0$ and propagate them in the 2D phase space. Second, calculate the joint and marginal density using the binning approach. For DEE, the main computational effort covers two aspects. The first is the generation of initial random samples at $t = 0$ and propagation of the samples in the 3D state space ($\phi$, $e$, $n$). Second, calculate the joint and marginal density with techniques of linear interpolation and binning approach. The main contribution of the computational effort of GMM-UT also includes two parts. First, split the initial Gaussian distribution into $N_{1d}$ sub-Gaussians, obtain initial $N_{1d} \cdot (2 \cdot Nvar + 1)$ sigma points of UT at $t = 0$ and propagate them in the 2D phase space. Second, according to the propagated sigma points, analytically calculate the combined joint and marginal density with the means and covariance matrices propagated using UT.

Based on the computational effort analysis, the main contribution of the computational effort is divided into two parts for each method. The first part is defined as 'Propagation' effort, referring to the propagation of dynamical equations. The second part is defined as '1/0 Interpolation', referring to the data processing and calculation for obtaining the joint and marginal density, and corresponding to the binning effort (0 Interpolation) of MC, or the interpolation and binning effort of DEE (1 Interpolation), or the computational effort of analytically calculating the combined joint and marginal density of GMM-UT (0 Interpolation). The computational effort analysis will be performed in terms of the computational time of the 'Propagation' part, $t_{prop}$, and the '1/0 Interpolation' part, $t_{1/0int}$, the total computational time of each case defined as the sum of the considered two-part computational effort, $t_{cal}$ ($t_{cal}=t_{prop}+t_{1/0int}$), the computational effort ratio of each method for the defined two-part computational effort, $Rt_{2p}$ ($Rt_{2p}=[t_{prop}/t_{cal}; t_{1/0int}/t_{cal}]$), and the normalised computational time of DEE and GMM-UT with respect to that of MC, $t_{cal}/t_{cal-MC}$, respectively. The simulation is done on a standard laptop with a 1.80 GHz Intel Core i7 processor and 8GB RAM.

## 5 Results and discussion

### 5.1 Scenario 1 within SubD1

In this section, the results in terms of the phase space evolution, joint density, marginal density and computational effort are presented for the test case scenarios presented in Tables 1 and 2. To assess the quality of the density propagation with DEE and GMM-UT, the relative errors of the predicted mean and standard deviation are also compared with MC. Note that detailed results and analyses are given for Scenario 1, while for Scenario 2 and Scenario 3, only the results in terms of the joint density, marginal density, relative errors of the mean and standard deviation are presented in "Appendices B and C", respectively, for parallel comparison with Scenario 1 to give an insight into density evolution results within each sub-phase space domain.

#### 5.1.1 Phase space evolution

Figure 2 shows the phase space evolution results for the first test case scenario with initial 961 and 1E5 Gaussian-distributed random samples depicted at times $t=\{0, 0.5, 1, 1.5, 2\}$ yrs. The two gray closed loop arcs denote the minimum (the inner one) and the maximum Hamiltonian contour line corresponding to the phase space random samples. It is shown that the phase space evolution results are within the sub-phase space domain of SubD$_1$ defined in Fig. 1. The results of $N_{sam}=961$ capture the characteristics of the phase space deformation and elongation with time well compared with that of the larger sample number of 1E5. Note that this is a significant requirement to ensure the performance of DEE with linear interpolation for density evolution, due to the dependence of the technique of linear interpolation on the quality of the samples to characterise the true distribution of phase space and density space. As seen in Fig. 2, for the case of an initial Gaussian density distribution under the influence of SRP and $J_2$, the deformation of phase space distribution increases with time and differs for different domains of the eccentricity. It is shown that at times $t=\{1, 2\}$ yrs, the phase space distribution appears much more elongated and deformed in the direction of the solar angle compared with that of the cases at times $t=\{0.5, 1.5\}$ yrs. This is mainly due to the higher nonlinearity of the dynamical system with respect to the specific eccentricity domain, either larger than the critical eccentricity $e_{cri}=0.6$ (shown in Fig. 1 and defined in Sect. 2.1) or approaching the zero eccentricity.

#### 5.1.2 Joint density

To compare the joint density evolution results for MC, DEE and GMM-UT, detailed results are shown together and separately at times $t=\{0, 0.5, 1, 1.5, 2\}$ yrs in Figs. 3 and 4, respectively, for the four cases of Scenario 1 (defined in Table 2). It is shown that the joint density evolution results are in the form of transformed non-Gaussian density distribution, but still featuring an obvious density core value like that of a Gaussian density distribution.

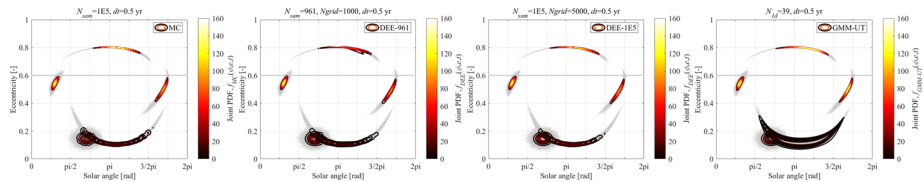

**Fig. 3** Integrated joint density for the four cases of Scenario 1 (critical eccentricity $e_{cri} = 0.6$ for Earth re-entry)

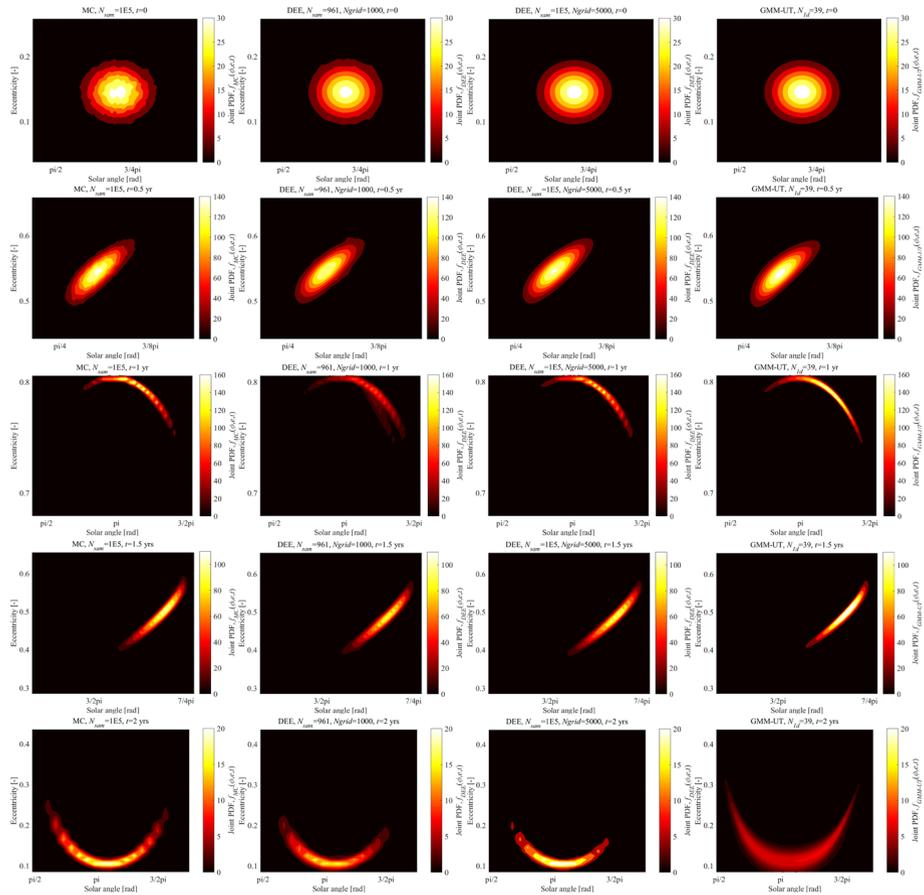

**Fig. 4** Separate joint density for the four cases of Scenario 1 at times $t = \{0, 0.5, 1, 1.5, 2\}$ yrs

For MC in Figs. 3 and 4, higher joint density peaks are obtained at times $t = \{0.5, 1, 1.5\}$ yrs than that at times $t = \{0, 2\}$ yrs. In Fig. 3, this can be seen by the colour contrast for each time instant, while in Fig. 4, this can be seen by the maximum colour data of the colourmap for each time instant. The maximum and minimum joint density peak are found at $t = 1$ yr with $(\phi, e) = (1.1618\pi, 0.7965)$ and at $t = 2$ yrs with $(\phi, e) = (1.1141\pi, 0.105)$, respectively.

For the accuracy of the joint density calculated by DEE and GMM-UT compared with MC, the following conclusions can be drawn from Figs. 3 and 4. First, for DEE-961, DEE-1E5 and GMM-UT (defined in Table 2), comparable and accurate joint density is obtained compared with that of MC in terms of the joint density distribution characteristics (in terms of the core value and the changing characteristics) shown in the integrated results for each time instant in Fig. 3, and the results presented separately in Fig. 4. Second, the joint density of DEE-1E5 is much more consistent with that of MC than that of DEE-961 and GMM-UT. This is mainly due to the same large number of samples as that of MC and a larger grid number $Ngrid$=5E3 to obtain linearly interpolated density. Third, for DEE-961, the joint density obtained using linear interpolation is highly consistent with that of MC, excluding the case of the underestimated result at $t$=1 yr. This is mainly due to the highly deformed and elongated phase space distribution that leads to worse interpolated density within the whole convex hull. Fourth, for GMM-UT, the joint density obtained is highly consistent with that of MC as shown in Figs. 3 and 4, except the overall underestimated result at $t$=2 yrs in terms of the core value. This is mainly due to the sensibility of GMM-UT to the highly nonlinear phase space dynamics when it approaches the zero eccentricity, and the operation in this paper of performing splitting only at the initial time $t$=0 along the direction of the solar angle.

It can be preliminarily concluded that for the case of less deformed and elongated phase space distribution, such as the cases of $t$<2 yrs in Figs. 3 and 4, GMM-UT outperforms DEE with linear interpolation in terms of the accuracy of the joint density. (Sect. 3.2 gives explanation of the technique of linear interpolation in terms of the calculation of density together with density evolution equation and binning approach techniques.) The performance of DEE is mainly constrained by the performance of the linear interpolation and the number of initial samples. The problem of obtaining better linearly interpolated results for the case of highly deformed and elongated phase space distribution (Trisolini and Colombo 2021) will be given into a deeper insight in the future work. For GMM-UT, the splitting of the initial Gaussian distribution is performed at the initial time along the solar angle direction. With time evolution, the splitted sub-Gaussians extend along the solar angle direction. This can be seen by referring to the last-column figures in Fig. 4 depicted at times $t$={0, 0.5, 1, 1.5, 2} yrs. Overall, there is an 80% probability that GMM-UT outperforms DEE in terms of capturing the core characteristic of the joint density. For example, at times $t$={0, 0.5, 1, 1.5} yrs, GMM-UT obtains a more continuous estimation of the joint density than that of MC. However, it fails to capture the elongated joint density characteristic when it approaches the singularity at zero eccentricity, such as the case of $t$=2 yrs.

### 5.1.3 Marginal density

To compare the marginal density evolution results of DEE and GMM-UT with MC, the results of DEE-961, DEE-1E5 and GMM-UT are shown in Fig. 5, together with that of MC, in terms of four comparison pairs of MC vs. DEE-961, MC vs. DEE-1E5, MC vs. GMM-UT and MC vs. DEE vs. GMM-UT.

For MC, the marginal density evolution histogram is represented in sequential blue colour and the marginal density evolution curve is given with a red line. It is shown that with time evolution, the marginal density distribution loses its initial Gaussian structure, especially for the cases at times $t = \{1, 2\}$ yrs for both solar angle and eccentricity. The maximum and minimum marginal density peak of the solar angle are obtained at $\phi = 0.3019\pi$ rad at $t = 0.5$ yr (marginal density value, $f_{MC\text{-}1\text{-}p} = 5.564$; unit of the marginal density of the solar angle: frequency of samples per bin width $wid_1$ (defined in Eq. (10))) and at $\phi = 0.7703\pi$ rad at $t = 2$ yrs ($f_{MC\text{-}1\text{-}p} = 0.4367$), respectively. Note that two local marginal density peaks exist at $t = 2$ yrs within the solar angle domain $[0.4832\pi, 1.7312\pi]$. Another local marginal density peak is at $\phi = 1.3207\pi$ rad ($f_{MC\text{-}1\text{-}p} = 0.3449$) at $t = 2$ yrs. Meanwhile, highly non-Gaussian characteristic appears at $t = 1$ yr in terms of the marginal density of the solar angle within the domain $[0.3988\pi, 1.5352\pi]$, where one peak is observed at $\phi = 1.2366\pi$ rad ($f_{MC\text{-}1\text{-}p} = 0.7436$). Overall, with the time evolution, the marginal density peak of the solar angle first increases until $t = 0.5$ yr and then drops after $t = 1$ yr. For the marginal density of the eccentricity, two maximum marginal density peaks are reached at $e = 0.8029$ at $t = 1$ yr (marginal density value, $f_{MC\text{-}2\text{-}p} = 101.6$; unit of the marginal density of the eccentricity: frequency of samples per bin width $wid_2$) and at $e = 0.1072$ at $t = 2$ yrs ($f_{MC\text{-}2\text{-}p} = 22.44$), respectively. This is consistent with the aforementioned higher nonlinearity characteristics of the dynamics when approaching an eccentricity value larger than the critical eccentricity $e_{cri}$ or equal to zero. The most non-Gaussian distribution characteristic is captured at $t = 1$ yr within the eccentricity domain $[0.6796, 0.8057]$. The minimum marginal density peak is obtained at $e = 0.5051$ at $t = 1.5$ yrs ($f_{MC\text{-}2\text{-}p} = 8.909$), indicating the minimum nonlinearity in the direction of the eccentricity at $t = 1.5$ yrs. Little discrepancy is shown for the marginal density distribution of the eccentricity at times $t = \{0, 0.5\}$ yr. Overall, with the time evolution, the marginal density peak of the eccentricity increases and the nonlinearity (or steepness) of the marginal density distribution increases when it approaches the eccentricity domain larger than the critical eccentricity or approaches the zero eccentricity.

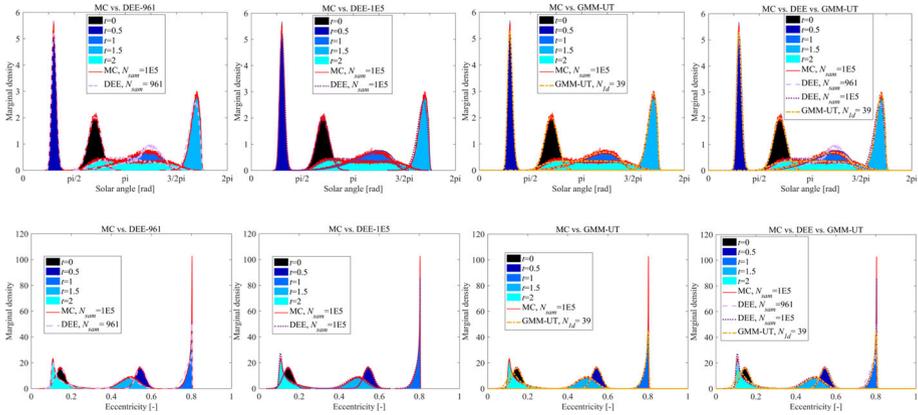

**Fig. 5** Marginal density of MC versus DEE-961, MC versus DEE-1E5, MC versus GMM-UT and MC versus DEE versus GMM-UT for Scenario 1

As shown in Fig. 5, the marginal density in terms of the solar angle and eccentricity for DEE-961 appears non-smooth at times $t = \{1, 2\}$ yrs. This is mainly due to the highly deformed phase space distribution, which complicates the linear interpolation of the density compared with the cases of the less deformed phase space at times $t = \{0, 0.5, 1.5\}$ yrs. Overall, the marginal density in terms of the solar angle and eccentricity is well captured compared with that of MC.

For DEE-1E5, the marginal density in terms of the solar angle and eccentricity shows a better consistency with respect to that of MC, since in this case the larger number of samples allows smoother marginal density. Note that the maximum value of the marginal density in terms of the eccentricity is overestimated compared with that of MC at $t=2$ yrs. As expected, a larger sample size generates a smoother and more consistent prediction of the uncertainty distribution as it improves the results of the linear interpolation. However, overestimation of the density peak may exist. This is due to the limitation of the linear interpolation to the convex hull of the scattered samples, which neglects parts of the distribution. Although these parts are associated with low density values, discarding them can lead to an overestimation of the peaks.

For the comparison between the MC method and the GMM-UT method, we can see that the marginal density obtained by GMM-UT in terms of both solar angle and eccentricity is well consistent with that of MC, except in correspondence of $t=2$ yrs. In this case, the marginal density of the eccentricity is underestimated, and the slope and skewness information are barely captured. Similar to the result of the underestimation of the joint density of GMM-UT at $t=2$ yrs shown in Fig. 4, this is mainly due to the higher nonlinearity of the phase space dynamics in the direction of the eccentricity when it approaches zero. Additionally, it has been considered that in this paper, the splitting of GMM is only performed at $t=0$ along the direction of the solar angle without considering splitting during the phase space propagation.

For the combined marginal density comparison results of MC vs. DEE vs. GMM-UT given in Fig. 5, it is shown that DEE-1E5 outperforms DEE-961 and GMM-UT in terms of the accuracy (including the slope and skewness information) of the marginal density of both solar angle and eccentricity with respect to that of MC. With a large number of initial samples of DEE (such as 1E5 here), the marginal density obtained using linear interpolation can be even smoother than the marginal density curve and density histogram obtained by MC.

### 5.1.4 Mean and standard deviation

To compare the accuracy of DEE and GMM-UT with MC for long-term density propagation in terms of the first two statistical moments, both the mean and standard deviation of the 2D phase space are given in Table 3 for the four cases of Scenario 1. The relative errors are depicted in Fig. 6 for DEE-961, DEE-1E5 and GMM-UT with respect to that of MC.

Table 3 Mean, standard deviation of the 2D phase space for the four cases of Scenario 1

|  | $t$, yr | $t=0$ | $t=0.5$ | $t=1$ | $t=1.5$ | $t=2$ |
|---|---|---|---|---|---|---|
| MC | $\mu_\phi$, rad | 2.2078 | 0.95229 | 3.57337 | 5.26674 | 3.32611 |
|  | $\sigma_\phi$, rad | 0.19671 | 0.07561 | 0.58152 | 0.14873 | 0.86184 |
|  | $\mu_e$ | 0.14499 | 0.54399 | 0.78955 | 0.49021 | 0.14285 |
|  | $\sigma_e$ | 0.025 | 0.0246 | 0.01501 | 0.04534 | 0.04342 |
| DEE-961 | $\mu_\phi$, rad | 2.20634 | 0.95963 | 3.66188 | 5.22771 | 3.23971 |
|  | $\sigma_\phi$, rad | 0.19407 | 0.07559 | 0.5835 | 0.15942 | 0.73567 |
|  | $\mu_e$ | 0.14596 | 0.54736 | 0.77975 | 0.47897 | 0.1488 |
|  | $\sigma_e$ | 0.02583 | 0.02428 | 0.01942 | 0.04554 | 0.05303 |
| DEE-1E5 | $\mu_\phi$, rad | 2.20692 | 0.95905 | 3.66475 | 5.23035 | 3.25211 |
|  | $\sigma_\phi$, rad | 0.1964 | 0.07676 | 0.56539 | 0.16053 | 0.75818 |
|  | $\mu_e$ | 0.14501 | 0.54713 | 0.78676 | 0.47854 | 0.13386 |
|  | $\sigma_e$ | 0.02511 | 0.02431 | 0.01715 | 0.04593 | 0.03863 |
| GMM-UT | $\mu_\phi$, rad | 2.20695 | 0.95208 | 3.57347 | 5.2672 | 3.33267 |
|  | $\sigma_\phi$, rad | 0.19635 | 0.07583 | 0.58318 | 0.15021 | 0.87551 |
|  | $\mu_e$ | 0.145 | 0.54399 | 0.78966 | 0.49027 | 0.1426 |
|  | $\sigma_e$ | 0.025 | 0.02467 | 0.01517 | 0.04551 | 0.04525 |

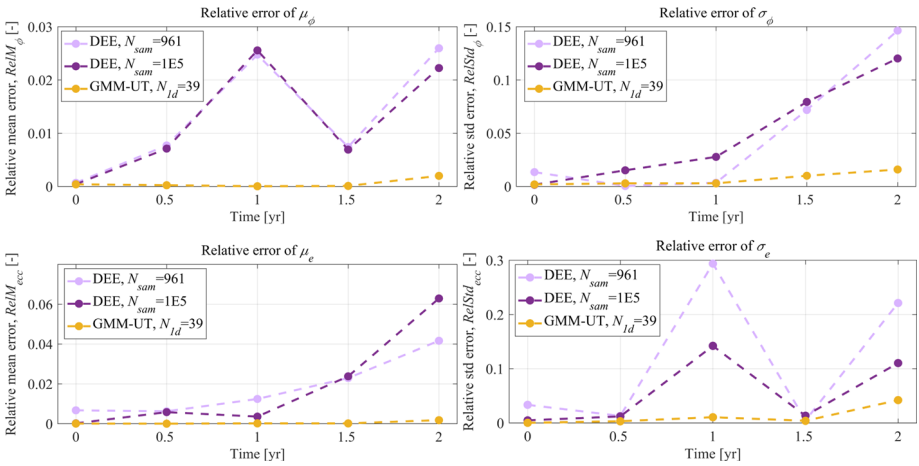

**Fig. 6** Relative error of mean, standard deviation of the solar angle; relative error of mean, standard deviation of the eccentricity (for Scenario 1)

As shown in Fig. 6, the relative errors in terms of the mean and standard deviation at different time instants are all under 7% and 30%, respectively. GMM-UT outperforms DEE with linear interpolation. Overall, DEE-1E5 obtains better consistency results than that of DEE-961 with respect to that of MC.

### 5.1.5 Computational effort

Figure 7 shows the computational effort ratios for the four cases of Scenario 1 for the defined two-part computational effort. Detailed results are given in Table 4 in terms of the parameter set of $\{t_{prop}, t_{1/0int}, t_{cal}, t_{cal}/t_{cal-MC}\}$ (defined in Sect. 4.2).

As shown in Table 4, the computational effort of DEE-1E5 is the heaviest for the four cases considered, over twice of that of MC, while that of GMM-UT is the lightest, accounting for only 0.42% of that of MC.

As depicted in Fig. 7, for MC, DEE-961 and DEE-1E5, the main contribution of the computational effort is the propagation of samples in the 2D ($\phi$, $e$) phase space and the 3D ($\phi$, $e$, $n$) state space, respectively, which is especially the case for MC and DEE-1E5. For MC, the minimum computational time of 0.41 s of the '1/0 Interpolation' part is taken to obtain the joint and marginal density since the calculation can be done with only binning approach. For DEE-1E5, the computational effort for the '1/0 Interpolation' part is the highest. In this case, both the binning approach and the linear interpolation method with $Ngrid = 5E3$ are applied. For GMM-UT, the main contribution of the computational effort is the '1/0 Interpolation' part. Also note that the computational effort for the '1/0 Interpolation' part of GMM-UT is less than that of DEE-961.

**Fig. 7** Computational effort ratios for the four cases of Scenario 1 for the defined two-part computational effort

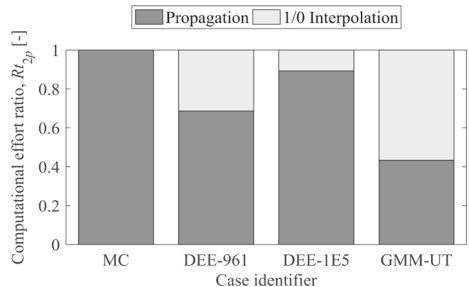

**Table 4** Results of computational effort for the four cases of Scenario 1

| Case | 1 | 2 | 3 | 4 |
|---|---|---|---|---|
| Identifier | MC | DEE-961 | DEE-1E5 | GMM-UT |
| $t_{prop}$, s | 4370.77 | 69.41 | 8712 | 7.95 |
| $t_{1/0int}$, s | 0.41 | 31.64 | 1048.83 | 10.39 |
| $t_{cal}$, s | 4371.17 | 101.05 | 9760.83 | 18.34 |
| $t_{cal}/(t_{cal-MC})$ | 1 | 0.02312 | 2.233 | 0.0042 |

### 5.2 Discussion

In this section, some discussion is made on the results in terms of the comparison of the density evolution results within sub-phase space domains and the comparison of MC, DEE and GMM-UT in terms of computational accuracy and efficiency.

#### 5.2.1 Density evolution result comparison within sub-phase space domains

Based on the density evolution results for the three test case scenarios (see Figs.3, 4,5 of Scenario 1 in Sect. 5.1, Figs. 10, 11,12 of Scenario 2 in "Appendix B" and Figs. 14, 15,16 of Scenario 3 in "Appendix C"), discussion is given on the comparison of the joint and marginal density evolution results within sub-phase space domains.

The following two-aspect results are drawn by comparing the density within each sub-phase space domain. First, considering the nonlinearity ranking of the phase space dynamics within each sub-phase space domain, i.e., {$SubD_1$, $SubD_2$}>$SubD_3$, a greater changing amplitude is obtained for Scenario 1 in Fig. 3 and Scenario 2 in Fig. 10 than that of Scenario 3 in Fig. 14 of the joint density peak shown in the integrated results. Also, weaker nonlinearity is detected for Scenario 3 in Fig. 16 and Scenario 2 in Fig. 12 compared with that of Scenario 1 in Fig. 5 of the marginal density evolution results for both solar angle and eccentricity. Second, note that the phase space distribution of Scenario 1 at $t$=1 yr in Fig. 2 and Scenario 2 at $t$=1.5 yrs in Fig. 10 are all highly deformed and elongated in the direction of the solar angle within the eccentricity domain larger than the critical eccentricity, but different characteristics of the joint density peak appear. For the former in Fig. 3, the maximum joint density peak is found, while for the latter in Fig. 10, the minimum joint density peak is found. This is mainly due to the fact that different dynamical characteristics are shown within each sub-phase space domain. It also points out the significance of a detailed analysis of the density evolution results within each sub-phase space domain.

### 5.2.2 Accuracy comparison of DEE and GMM-UT

As shown in the detailed accuracy analysis of Scenario 1 in Sect. 5.1 and results of Scenarios 2-3 in "Appendices B, C" in terms of the joint density, marginal density, mean and standard deviation for both DEE and GMM-UT compared with that of MC (see Figs. 3, 4, 5,6 of Scenario 1, Figs.10, 11, 12,13 of Scenario 2 and Figs.14, 15, 16,17 of Scenario 3), overall, both methods work well for long-term density propagation within each sub-phase space domain.

For the joint or marginal density evolution, for cases where dynamical nonlinearity is high, such as Scenario 1 within $SubD_1$, the DEE with linear interpolation using a large number of initial samples $N_{sam}$=1E5 outperforms GMM-UT and obtains highly consistent density with respect to that of MC (but with a much higher computational load compared with GMM-UT). Note that the slightly worse performance of GMM-UT here for the case exhibiting higher dynamical nonlinearity is probably due to the operation in this paper of performing Gaussian splitting only at the initial time along the solar angle direction. For cases where the dynamical nonlinearity is lower, such as Scenario 2 and Scenario 3, the performance of GMM-UT is comparable to that of DEE-1E5 (but with a much lower computational load). For the accuracy of the propagation of the first two statistical moments, GMM-UT outperforms DEE.

The accuracy of DEE with linear interpolation mainly depends on the dynamical nonlinearity and the performance of linear interpolation. To improve the performance of linear interpolation when the phase space is highly deformed and elongated, an improved linear interpolation method considering alpha shape (Trisolini and Colombo 2021) can be studied in the future work. The accuracy of GMM-UT mainly depends on the dynamical nonlinearity, the splitting number of sub-Gaussians (i.e., theoretically, as the number of splitting increases, a GMM can approach the true non-Gaussian density distribution), the splitting direction (i.e., whether the splitting is done along the most nonlinear direction of the dynamical system), whether the splitting is done within various directions at the initial time (Vittaldev et al. 2016), whether the splitting is done during the propagation by detecting the nonlinearity of the state uncertainty (DeMars et al. 2013; Romano 2020), and the performance of UT for nonlinear propagation of the first two statistical moments of each sub-Gaussian. These aspects need to be given into an insight in the future work to improve the performance of GMM-UT.

### 5.2.3 Efficiency comparison of DEE and GMM-UT

Figure 8 shows a comparison of the normalised computational effort for the three test case scenarios. It is shown that for all sub-phase space domains, the computational effort of DEE-1E5 is the heaviest, which is often over twice of that of MC (except the case of Scenario 2), while that of GMM-UT is the lightest. For DEE-1E5, this is due to the fact that the same large number as that of MC of initial samples is propagated in the 3D state space ($\phi$, $e$, $n$). The linear interpolation method with *Ngrid* = 5E3 is applied to calculate the joint and marginal density. For GMM-UT, this is due to the facts that only 195 ($N_{1d}\cdot(2\cdot Nvar + 1)$, $N_{1d} = 39$, $Nvar = 2$) sigma points are propagated, the mean and covariance matrix corresponding to each sub-Gaussian are calculated with explicit equations of unscented transformation, and the joint and marginal density are analytically calculated using the weighted results of the mean and covariance matrix of each sub-Gaussian. The computational effort of DEE-961 is 3.5 to 5 times of that of GMM-UT. This is mainly due to the computational effort of propagating nearly five times the number of initial samples of that of GMM-UT in the 3D state space, and the application of the linear interpolation method with *Ngrid* = 1E3 to obtain the joint and marginal density.

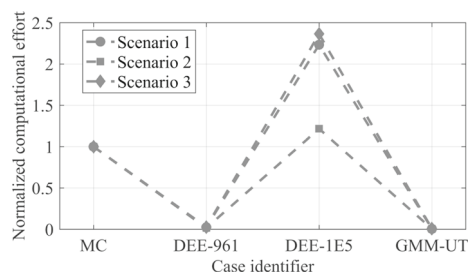

**Fig. 8** Comparison of normalised computational effort for the four cases of the three test case scenarios

## 6 Conclusion

This paper compares Density Evolution Equation (DEE) and Gaussian Mixture Model (GMM) on the long-term 2D phase space density propagation problem in the context of high-altitude and high area-to-mass ratio satellite long-term propagation. The density evolution equation is formulated under the influence of solar radiation pressure and Earth's oblateness using the semi-analytical method. The linear interpolation method based on Delaunay triangulation is integrated with density evolution equation for accurate density calculation. Unscented Transformation (UT) is used to nonlinearly propagate the first two statistical moments corresponding to each sub-Gaussian. An insight is given into the 2D phase space long-term density propagation problem within three sub-phase space domains subject to non-linear dynamics. For the propagation accuracy of the joint and marginal density, overall, both DEE and GMM-UT work well within each sub-phase space domain. For cases with higher dynamical nonlinearity, the DEE with linear interpolation using a great number of initial samples $N_{sam}$=1E5 outperforms GMM-UT and obtains highly consistent density with respect to that of MC (but with a much higher computational load compared with GMM-UT). For cases with lower dynamical nonlinearity, the performance of GMM-UT is comparable to that of DEE (but with a much lower computational load). For the propagation accuracy of the first two statistical moments, GMM-UT outperforms DEE. Overall, GMM-UT is more suitable for the long-term density propagation problem. Future work on improving the performance of GMM-UT will provide a better solution for the long-term density propagation problem.

**Acknowledgements** This project has received funding from the European Research Council under the European Union's Horizon 2020 research and innovation programme (grant agreement 679086–COMPASS). Dr. Pan Sun and Shuang Li thank the National Natural Science Foundation of China (Grant No. 11972182 and 11672126), Qing Lan Project, China Scholarship Council (CSC No. 202006830123), and the 2020 Postgraduate Research Practice Innovation Program of Jiangsu Province (Grant No. KYCX20_0222), which funded their participation to this research.

## Appendix A: Geometrical representation of the problem

Four orbital elements (semi-major axis $a$, eccentricity $e$, longitude of pericentre $\tilde{\omega}$ and mean anomaly $M$) are required to describe the planar problem. In this paper, the shadowing effects are not considered, allowing to treat the semi-major axis $a$ as constant (Allan 1962). Similar to the works of (Krivov and Getino 1997) and (Lücking et al. 2011a), we replace the longitude of pericentre $\tilde{\omega}$ with the more physically meaningful variable of the solar angle $\phi$. The two variables are connected via the equation $\phi = \tilde{\omega} - \lambda_{sun}$, where $\lambda_{sun}$ is the true longitude of the sun. It represents the angular distance between the pericentre and the direction toward the sun. In this paper, zero Earth orbital eccentricity is assumed, and thus, $\lambda_{sun}$ is a linear function of time $t$, i.e. $\lambda_{sun} = n_{sun} \cdot t$, where $n_{sun}$ is the mean motion of the sun. Figure 9 shows a geometrical representation of the problem, where $\Gamma$ direction targets at the vernal equinox.

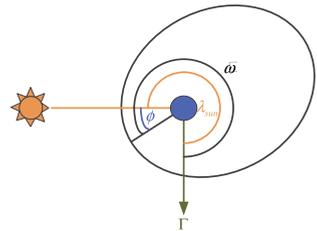

**Fig. 9** Planar Earth orbit in rotating reference frame (Lücking et al. 2011a)

## Appendix B: Scenario 2 within SubD2

Figure 10 shows the joint density together at times $t=\{0, 0.5, 1, 1.5, 2, 2.5, 3\}$ yrs for the four cases of Scenario 2 within sub-phase space domain SubD$_2$. Figure 11 shows the joint density separately at times $t=\{0.5, 1, 1.5, 2, 2.5, 3\}$ yrs. Detailed marginal density evolution results of DEE-961, DEE-1E5 and GMM-UT are shown in Fig. 12, together with that of MC, in terms of four comparison pairs of MC vs. DEE-961, MC vs. DEE-1E5, MC vs. GMM-UT and MC vs. DEE vs. GMM-UT. The relative errors in terms of the mean and standard deviation for DEE and GMM-UT with respect to MC are depicted in Fig. 13.

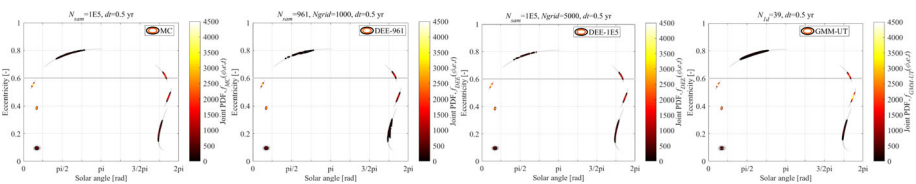

**Fig. 10** Integrated joint density for the four cases of Scenario 2 (critical eccentricity ecri = 0.6 for Earth re-entry)

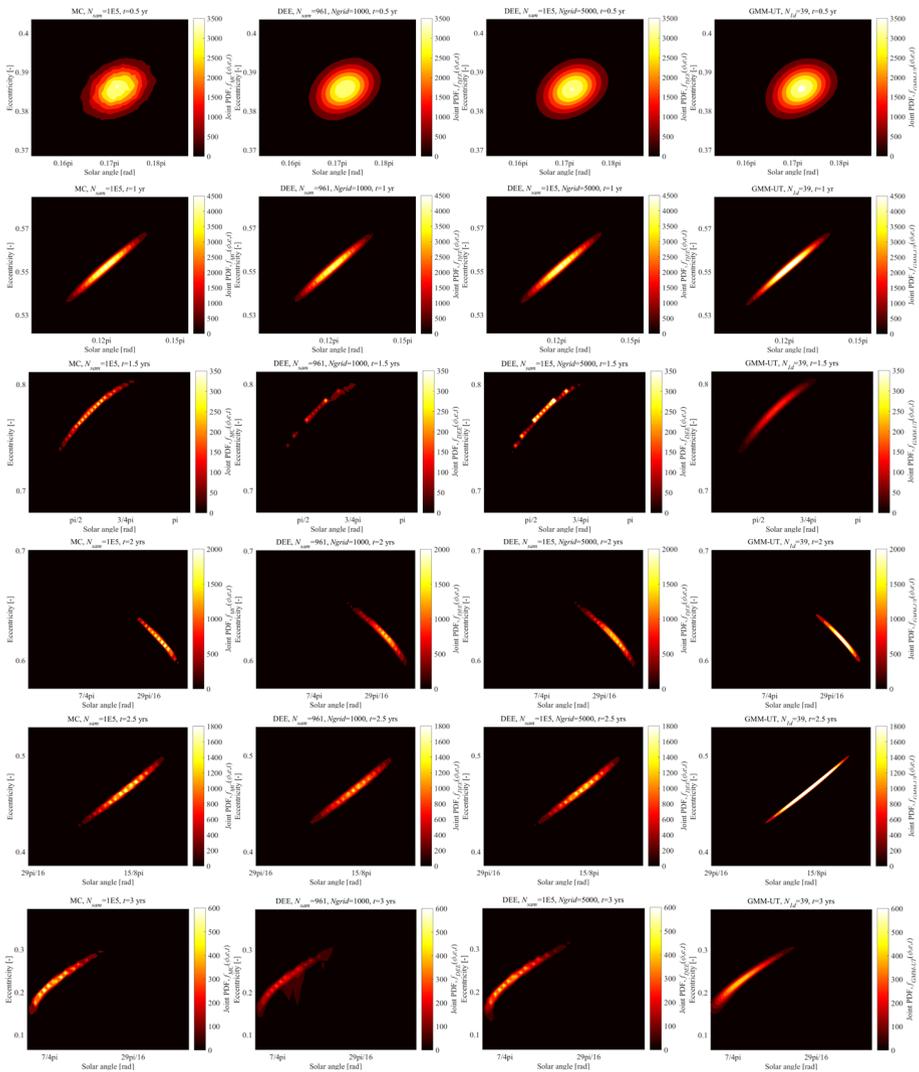

**Fig. 11** Separate joint density for the four cases of Scenario 2 at times $t = \{0.5, 1, 1.5, 2, 2.5, 3\}$ yrs

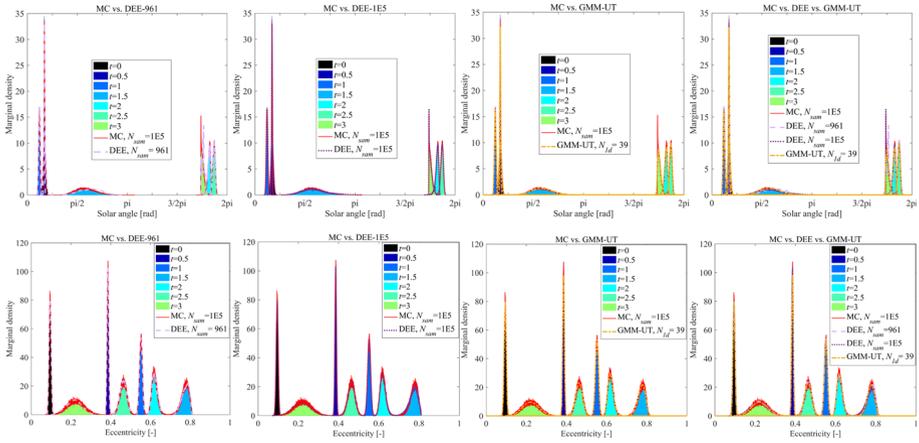

**Fig. 12** Marginal density of MC versus DEE-961, MC versus DEE-1E5, MC versus GMM-UT and MC versus DEE versus GMM-UT for Scenario 2

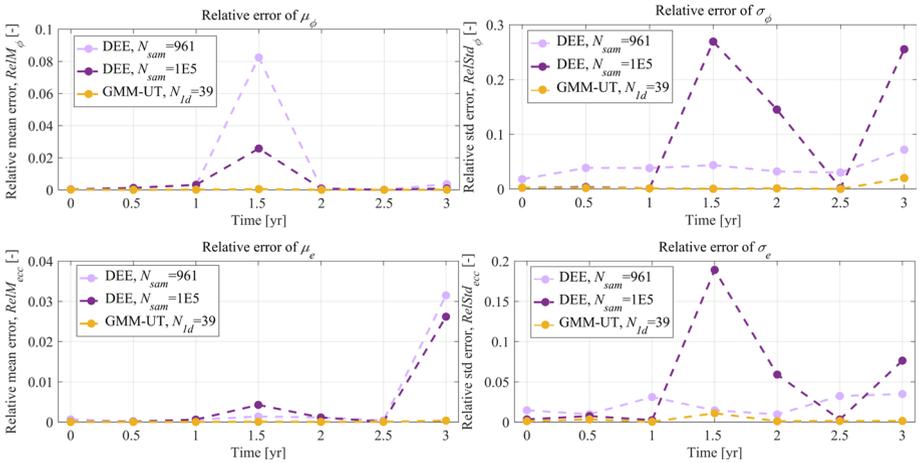

**Fig. 13** Relative error of mean, standard deviation of the solar angle; relative error of mean, standard deviation of the eccentricity (for Scenario 2)

## Appendix C: Scenario 3 within SubD3

Figure 14 shows the joint density together at times $t=\{0, 0.5, 1, 1.5, 2\}$ yrs for the four cases of Scenario 3 within sub-phase space domain SubD$_3$. Figure 15 shows the joint density separately at times $t=\{0.5, 1, 1.5, 2\}$ yrs. Here the joint density is presented within the solar angle domain $[-\pi, \pi)$. Detailed marginal density evolution results of DEE-961, DEE-1E5 and

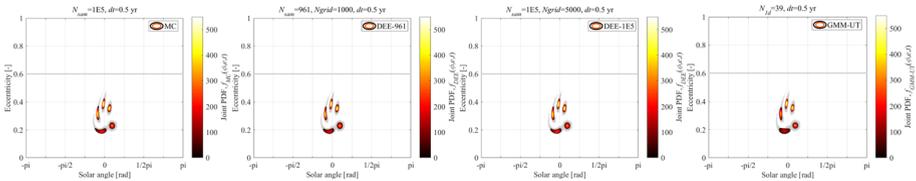

**Fig. 14** Integrated joint density for the four cases of Scenario 3 (critical eccentricity $e_{cri} = 0.6$ for Earth re-entry)

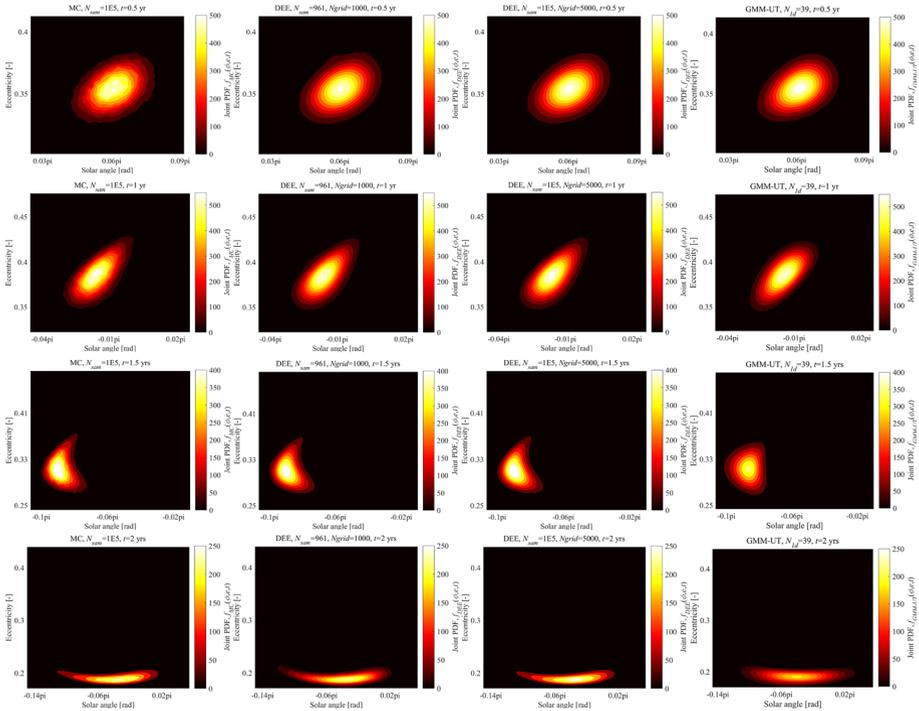

**Fig. 15** Separate joint density for the four cases of Scenario 3 at times $t = \{0.5, 1, 1.5, 2\}$ yrs

GMM-UT are shown in Fig. 16, together with that of MC, in terms of four comparison pairs of MC vs. DEE-961, MC vs. DEE-1E5, MC vs. GMM-UT and MC vs. DEE vs. GMM-UT. The relative errors in terms of the mean and standard deviation for DEE and GMM-UT with respect to MC are depicted in Fig. 17.

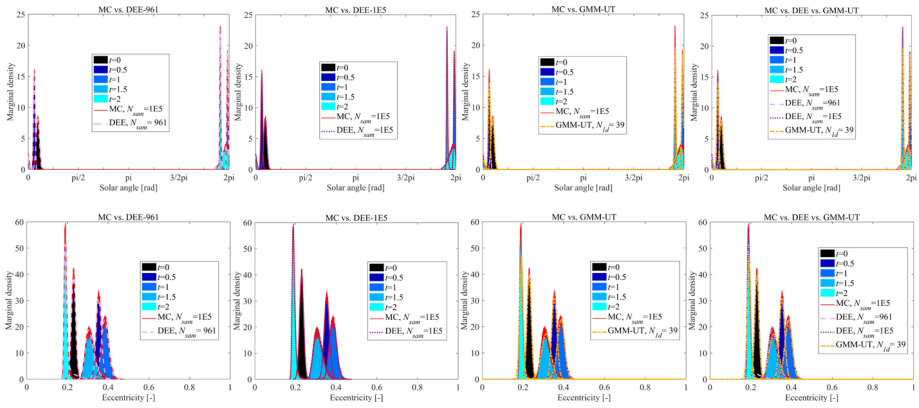

**Fig. 16** Marginal density of MC versus DEE-961, MC versus DEE-1E5, MC versus GMM-UT and MC versus DEE versus GMM-UT for Scenario 3

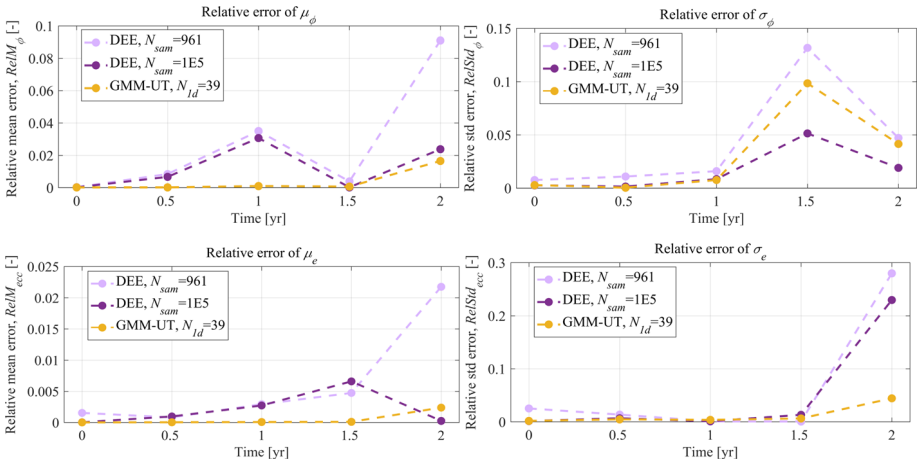

**Fig. 17** Relative error of mean, standard deviation of the solar angle; relative error of mean, standard deviation of the eccentricity (for Scenario 3)